\documentclass[hyperref=true]{elsarticle}
\usepackage[english]{babel}
\usepackage[utf8]{inputenc}
\usepackage[margin=2.5cm]{geometry}
\usepackage{natbib}
\usepackage[backref]{hyperref}
\hypersetup{colorlinks, citecolor=black, linkcolor=black, urlcolor=black, backref}
\usepackage{mathtools} 

\usepackage{cleveref}

\usepackage{float}

\usepackage{amssymb}  
\usepackage{amsthm}
\usepackage{amsmath}
\usepackage{commath}
\usepackage{booktabs}
\usepackage{placeins} 
\usepackage{subcaption}
\usepackage{framed}
\usepackage[normalem]{ulem}
\usepackage{enumerate}
\usepackage{cases}

\usepackage{todonotes} 
\newtheorem{lemma}{Lemma}

\usepackage{tikz-cd}

\newcommand{\Id}{I}
\newcommand{\uu}{{u}}
\renewcommand{\AA}{\mathcal{A}}
\newcommand{\Aapp}{\mathbb{A}}

\newcommand{\LL}{\mathcal{L}}
\newcommand{\Ai}{\mathcal{A}_i}
\newcommand{\Ae}{\mathcal{A}_e}

\newcommand{\vv}{{v}}

\newcommand{\ff}{{f}}

\newcommand{\divg}{\nabla\cdot}

\newcommand{\diverg}{\nabla\cdot}

\theoremstyle{definition}
\newtheorem{theorem}{Theorem}
\newtheorem{remark}{Remark}

\graphicspath{{./}{./figures/}}

\begin{document}

\begin{frontmatter}




\title{Splitting Schemes for Coupled Differential Equations\\
{\small Block Schur-Based Approaches \& Partial Jacobi Approximation}
} 


\author[kaust]{Roberto Nuca}
\author[bergen,forde]{Erlend Storvik}
\author[bergen]{Florin A. Radu }
\author[nottingham]{Matteo Icardi\corref{cor1}}
\cortext[cor1]{matteo.icardi@nottingham.ac.uk}

\address[kaust]{CEMSE, King Abdullah University of Science and Technology, Thuwal, Saudi Arabia}
\address[bergen]{Center for Modeling of Coupled Subsurface Dynamics, Department of Mathematics,
University of Bergen, All\'egaten 44, 5007 Bergen, Norway}
\address[forde]{Department of Computer Science, Electrical Engineering and Mathematical Sciences, Western Norway University of Applied Sciences, Svanehaugvegen 1, 6812 Førde, Norway}
\address[nottingham]{School of Mathematical Sciences, University of Nottingham, NG7 2RD, Nottingham, UK}

\begin{abstract}
    Coupled multi-physics problems are encountered in countless applications and pose significant numerical challenges. In a broad sense, one can categorise the numerical solution strategies for coupled problems into two classes: monolithic approaches and sequential (also known as split, decoupled, partitioned or segregated) approaches. The monolithic approaches treat the entire problem as one, whereas the sequential approaches are iterative decoupling techniques where the different sub-problems are treated separately. Although the monolithic approaches often offer the most robust solution strategies, they tend to require ad-hoc preconditioners and numerical implementations. Sequential methods, on the other hand, offer the possibility to add and remove equations from the model flexibly and rely on existing black-box solvers for each specific equation. Furthermore, when problems are non-linear, inner iterations need to be performed even in monolithic solvers, making the sequential approaches an even more viable alternative. The cost of running inner iterations to recover the multi-physics coupling could, however, easily become prohibitive. Moreover, the sequential approaches might not converge at all. 
    In this work, we present a general formulation of splitting schemes for continuous operators with arbitrary implicit/explicit splitting, like in standard iterative methods for linear systems. By introducing a generic relaxation operator, we find the conditions for the convergence of the iterative schemes. We show how the relaxation operator can be thought of as a preconditioner and constructed based on an approximate Schur complement. We propose a Schur-based Partial Jacobi relaxation operator to stabilise the coupling and show its effectiveness. Although we mainly focus on scalar-scalar linear problems, most results are easily extended to non-linear and higher-dimensional problems. The schemes presented are not explicitly dependent on any particular discretisation methodologies. Numerical tests (1D and 2D) for two PDE systems, namely the Dual-Porosity model and a Quad-Laplacian operator, are carried out to investigate the practical implications of the theoretical results.
\end{abstract}

\begin{keyword}
partitioned multi-physics, splitting schemes, approximate Schur complement, partial Jacobi, block-iterative schemes, sequential coupling
\end{keyword}
\end{frontmatter}

\section{Introduction} \label{sec:intro}
Multi-physics phenomena can involve large systems of linear or non-linear partial differential equations. Important examples of coupled multiphysics problems include bulk coupling (where equations are solved in the same domain), such as Navier-Stokes equations, poromechanics, multi-species reactive transport, and interface couplings such as fluid-structure interaction and conjugate heat transfer, to name a few. Many theoretical and numerical studies have been devoted to these problems, resulting in optimised monolithic solution approaches that often need detailed analysis of the system, ad-hoc coupled stable discretisation techniques and pre-conditioners for the resulting (large) linear systems. This is not always achievable in computational engineering, where new coupled models often arise from the need to couple existing sub-models as black-boxes. In this direction, there have recently been efforts in developing generic frameworks capable of combining existing solvers, see \citep{preciceBUNGARTZ2016250}. These provide a layer of abstraction where solvers communicate with proper data mapping.
    
The strategies adopted to solve linear systems derived by the discretisation of coupled problems can be divided into monolithic and segregated/partitioned approaches. In monolithic approaches, the whole system is discretised in a single matrix, and all unknowns are solved simultaneously. In contrast, sequential iterative approaches consider each equation/sub-problem sequentially in a sequence with an internal loop to recover the coupled structure of the full problem. We denote here by "splitting scheme", the iterative and partitioned strategy that defines this internal loop. Internal iterations in splitting schemes can be conveniently combined with the non-linear iterations (Picard or Newton) required to solve non-linear problems.

Besides pure numerical aspects, there are software engineering motivations for why segregated/partitioned approaches could be preferable to monolithic approaches. One motivation is related to memory management because only one matrix at a time, derived by one single equation, needs to be stored. A second and probably most relevant motivation is the possibility of using many efficient and robust solvers already available in the scientific computing community to solve each equation appearing in the multiphysics model.

In this work, we focus on coupled scalar-scalar linear model problems, although most of the results are valid or easily extendable to other linear problems. The objective of this work is to develop a new family of splitting schemes based on approximate Schur complements and provide a theoretical framework to analyse them.
The new  framework generalises standard stabilised decoupling methods such as the fixed-stress 
splitting scheme \cite{settari1998, kim2011fixedstress, mikelicwheeler, jakubaml} and the undrained splitting scheme \cite{kim2011undrained} for the Biot 
equations of poroelasticity. These are based on adding a diagonal relaxation operator to stabilise the sequential method.
Recent works have  discussed the optimisation of this relaxation term for the fixed-stress 
splitting scheme in \cite{storvik2019optimization}.
Here, one of the overarching goals is to provide some insight into how more general operators
can be used to stabilise coupled operators.
While only linear problems are studied here, the generalisation 
to non-linear partial differential equations is straightforward as the splitting approaches proposed here can be easily
combined with linearisation methods, such as those applied for Richards equation coupled with transport in \cite{illiano2021iterative}.
Moreover, as with all fixed-point iterations, the methods proposed can be further improved by acceleration methods, such 
as the Anderson acceleration, which is known to accelerate linearly convergent schemes \cite{evans2020proof}. This has already
been successfully applied in combination with decoupling methods for the Biot equations \cite{both2019anderson} and phase-field methods 
for brittle fracture propagation \cite{storvik2021accelerated}.

The paper is organised as follows. In \cref{sec:model}, we introduce two model problems, the dual-porosity problem and the Quad-Laplacian. In the next section, we present a unified nomenclature for decoupling methods, conveniently seen as block-based versions of standard iterative schemes for linear systems.  We split the monolithic problem into an implicit and an explicit operator and then introduce a relaxation operator to improve and stabilise the iterations. In \cref{sec:schur}, we extend the idea behind the Uzawa algorithm and propose a new general decoupling method for solving iteratively coupled systems of equations based on an approximation of the Schur complement. In \cref{sec:convergence}, convergence conditions for generic splitting schemes and convergence proofs for several special cases, namely the block-symmetric and block-skew-symmetric cases, are provided for the new method. Eventually, in \cref{sec:results}, numerical results with the Schur-based Partial Jacobi approximation are compared to standard iterative methods.


\section{Model problems}\label{sec:model}
In this section, we present two model problems as examples of a system of partial differential equations falling in the generic representation used in \cref{sec:iterative_methods}, showing the conditions for well-posedness. In the last chapter, we will use these models as examples in numerical tests. 

\subsection{Dual-porosity Darcy's flow}
\label{sec:DualDarcy}
The first model problem we introduce is the so-called "dual-porosity" Darcy's flow model \cite{douglas1990dual}, which is an effective/mixture model to describe
the flow through a fractured porous medium.
Here, $u$ and $v$ are the matrix and fracture pressures, respectively. The equations read as follows:
\begin{align}
-\diverg\left(m_u \nabla u\right) = \beta(v-u) + \ff_1 \,, \label{eq:darcy1}\\
-\diverg\left(m_v\nabla v\right) = \beta(u-v) + \ff_2 \,, \nonumber
\end{align}
where $m_{(\cdot)}$ is the \textit{mobility} (permeability divided by the fluid viscosity), $\beta$ is the \textit{transfer coefficient} and $\ff_{1,2}$ are source terms. For the sake of simplicity, we consider homogeneous Dirichlet boundary conditions, but the methods developed in this paper can be applied to more general boundary conditions as well.

The continuous variational formulation of problem \cref{eq:darcy1} reads:\\

\noindent Find $(u,v) \in H^1_0(\Omega) \times H^1_0(\Omega)$ such that
\begin{align}
    a(u,\varphi) + b(u,\varphi) - b(v,\varphi) = \langle \ff_1,\varphi\rangle\qquad\forall\varphi\in H^1_0(\Omega)\,, \label{eq:darcy1_weak}\\
    d(v,\vartheta) + b(v,\vartheta) - b(u,\vartheta) = \langle \ff_2,\vartheta\rangle\qquad\forall\vartheta\in H^1_0(\Omega)\,,  \nonumber
\end{align}
where
$$
    a(u,\varphi) := \int_\Omega m_u\nabla u\cdot\nabla\varphi\,dx
    \qquad
    d(v,\vartheta) := \int_\Omega m_v\nabla v\cdot\nabla\vartheta\,dx
    \qquad
    b(u,\varphi) := \int_\Omega \beta u\varphi\,dx.
$$
We can sum \cref{eq:darcy1_weak} and using the bi-linearity property of the form $b(\cdot,\cdot)$ to get
$$
    a(u,\varphi) + b(u-v,\varphi-\vartheta) + d(v,\vartheta) = \langle \ff_1,\varphi\rangle + \langle \ff_2,\vartheta\rangle\qquad\forall(\varphi,\vartheta)\in H^1_0(\Omega)\times H^1_0(\Omega).
$$
The form 
$$
    \mathcal{A}:W\times W\to\mathbb{R}\qquad W:=H^1_0(\Omega)\times H^1_0(\Omega)
$$
$$
    ((u,v),(\varphi,\vartheta)) \mapsto a(u,\varphi) + d(v,\vartheta) + b(u-v,\varphi-\vartheta)
$$
is continuous and coercive under the following conditions:
$$
    m_u,m_v>0\,,\beta\geq0\text{ a.e. and }m_u,m_v\in L^\infty(\Omega).
$$
Hence, the problem \cref{eq:darcy1} is well-posed thanks to the Lax-Milgram Lemma. We will use this problem as a prototype of a block-symmetric coupled operator.
\subsection{Quad-Laplacian}
\label{sec:quadLaplacian}
The second problem we study is composed of four Laplace operators, and we name it as Quad-Laplacian. The main motivation for introducing this synthetic model is to study scalar problems with non-trivial off-diagonal "coupling" operators. 
The strong form of Quad-Laplacian reads
\begin{align}
    -\diverg\left(m_{uu} \nabla u\right) -\diverg\left(m_{uv} \nabla v\right) = f_1\,,\label{eq:quadLaplacian1}\\
    -\diverg\left(m_{vu}\nabla u\right) -\diverg\left(m_{vv} \nabla v\right) = f_2\,, \nonumber
\end{align}
where $m_{(\cdot,\cdot)}$ are diffusivity coefficients. Its weak form we get
\begin{align}
    a(u,\varphi) + b(v,\varphi) = \langle \ff_1,\varphi\rangle\qquad\forall\varphi\in H^1_0(\Omega)\,, \label{eq:quadLaplacian1_weak}\\
    c(u,\vartheta) + d(v,\vartheta) = \langle \ff_2,\vartheta\rangle\qquad\forall\vartheta\in H^1_0(\Omega)\,, \nonumber
\end{align}
$$
a(u,\varphi) = \int_\Omega m_{uu}\nabla u\cdot\nabla\varphi\,dx,
\qquad
b(v,\varphi) = \int_\Omega m_{uv}\nabla v\cdot\nabla\varphi\,dx,
$$
$$
c(u,\vartheta) = \int_\Omega m_{vu}\nabla u\cdot\nabla\vartheta\,dx,
\qquad
d(v,\vartheta) = \int_\Omega m_{vv}\nabla v\cdot\nabla\vartheta\,dx.
$$
Assuming $m_{uv}=-m_{vu}$, i.e., $b(v,u)=-c(u,v)$, 
$$\mathcal{A}:W\times W\to\mathbb{R}\qquad W:=H^1_0(\Omega)\times H^1_0(\Omega)$$
$$((u,v),(\varphi,\vartheta)) \mapsto a(u,\varphi) + d(v,\vartheta) + b(v,\varphi) + c(u,\vartheta)$$
is continuous and coercive under the following conditions:
$$m_{uu},m_{vv}>0\text{ a.e. and }m_{uu},m_{vv},m_{uv}\in L^\infty(\Omega).$$
Hence, the problem \cref{eq:quadLaplacian1} is well-posed due to the Lax-Milgram Lemma.
In the following, we will assume this condition to be valid, making the system block-skew-symmetric. 
This structure is typical in many coupled problems arising from mass and momentum conservation and semi-discrete poroelasticity.
For the quad-Laplacian case, each operator
is also symmetric (i.e., self-adjoint), making the whole system skew-symmetric.

\subsection{Well-posedness conditions and connection with saddle-point problems}

The models presented above are special cases of the following system of equations: 
\begin{equation}
\label{eq:system}
    \mathcal{A}w=
    \begin{bmatrix}
    	A & B \\
    	C & D 
    \end{bmatrix}
   	\begin{bmatrix}
    	\uu\\
    	\vv
	\end{bmatrix}
    =
    \begin{bmatrix}
    	\ff_1\\
    	\ff_2
	\end{bmatrix}
	=\ff 
\end{equation}
where $A$, $B$, $C$ and $D$ represent linear differential operators in the strong form.
This system is intentionally left without formal definitions of operators and involved spaces, as such details are discussed in the next section. In this work, we mostly adopt the terminology and concepts from linear algebra theory. Since not all steps are equally applicable to infinite dimensional (e.g. differential) operators, we occasionally interpret and indicate specific assumptions and differences for the case of differential operators separately. System \cref{eq:system} for Hilbert spaces is very general and a  simple approach to studying the well-posedness conditions, as done in the previous section for the Dual-Porosity and Quad-Laplacian system, is deriving a bilinear form acting between the Cartesian product of proper spaces and applying the hypotheses of Lax-Milgram lemma. However, one of the most relevant cases of \cref{eq:system} is the case of $D=0$ and $B=C^\top$ corresponding to the well-known saddle point problem (SPP) and its generalisations. For SPP, a robust and comprehensive theory has been developed \citep{fortin1991mixed,boffi2013mixed}. A generalisation of SPP where $D\neq0$ is sometimes referred to as a perturbed saddle point problem like those arising for nearly incompressible materials, recently analysed by \citep{hong2021new}. The SPP and perturbed SPP are examples where applying the Lax-Milgram lemma provides an unsatisfactory stability estimate because of the nature of the perturbative term \citep[Sec.~4.3, p.~238]{boffi2013mixed}. Another generalisation of SPP is where $D=0$ but $B$ and $C$ are generic. This case has been studied in \citep{nicolaides1982existence} as a direct generalisation of Brezzi-Babuska Theory. Even non-variational methods can be used to study the well-posedness of \cref{eq:system} provided proper spaces for the monolithic operator $\mathcal{A}$ \citep[Ch.~5]{atkinson2009theoretical}. In the context of PDEs, the continuity of the operators is a natural hypothesis, while coercivity and strong monotonicity for variational and non-variational methods, respectively, is the key hypothesis to ensure the injectivity of the operator.\\

In the continuous setting, \cref{eq:system} makes sense once proper spaces have been defined. Let $U$ and $V$ be Hilbert spaces with $u\in U$ and $v\in V$. We have that
\[
A:U\longrightarrow U'\,,\qquad B:V\longrightarrow U'
\,,\qquad 
C:U\longrightarrow V'\,,\qquad D:V\longrightarrow V'\,,
\]
hence $\ff_1\in U'$ and $\ff_2\in V'$. We can use the duality pairing in both equations and sum:
\[
\langle Au,\varphi\rangle_\ast + \langle Bv,\varphi\rangle_\ast + \langle Cu,\mu\rangle_\star + \langle Dv,\mu\rangle_\star = \langle \ff_1,\varphi\rangle_\ast + \langle \ff_2,\mu\rangle_\star\qquad\forall(\varphi,\mu)\in U\times V,
\]
i.e.,
\[
a(u,\varphi) + b(v,\varphi) + c(u,\mu) + d(v,\mu) = \langle \ff_1,\varphi\rangle_\star + \langle \ff_2,\mu\rangle_\star\qquad\forall(\varphi,\mu)\in U\times V.
\]

This expression defines a bilinear form
$$\mathcal{A}:W\times W\longrightarrow\mathbb{R}$$
\[
\mathcal{A}(w,\eta)=\langle\ff,\eta\rangle\qquad\forall \eta\in U\times V\,,
\]
where $W:=U\times V$ and $\ff$ is the functional defined by the sum of $\ff_1$ and $\ff_2$.
In the following, we will always assume $a(\cdot,\cdot)$ and $d(\cdot,\cdot)$ to be continuous and coercive with coercivity constants $\alpha_a$ and $\alpha_d$, respectively (i.e., the discrete operator $A$ and $D$ to be positive definite), and the monolithic operator
$\mathcal{A}$ to be trivially well-posed by assuming its continuity and coercivity. This can be a consequence of the hypotheses introduced for the dual-porosity or the quad-Laplacian model problems, or, more in general, by assuming that there exists a  coercivity constant $\alpha_\AA>0$ such that
\begin{equation}
\label{eq:monolithic_coerc_cond}
 0<\alpha_\AA\leq \mbox{min}(\alpha_a,\alpha_d)-\frac{\|b\|+\|c\|}{2}\,,
\end{equation}
where $\alpha_a$ and $\alpha_d$ are the coercivity constants of the bilinear forms $a(\cdot,\cdot)$ and $d(\cdot,\cdot)$ respectively.
In fact we have that, $\forall u,v \in U\times V$,
%
%
%
%
%

\[
\left| b(v,u) + c(u,v)\right| \leq \left(\|b\|+\|c\|\right)\|v\|\|u\|\leq \frac{\|b\|+\|c\|}{2}\left(\|v\|^2 + \|u\|^2\right).
 \]
Throughout the paper, although something else is specified, we use Euclidean norms for vectors and the induced matrix norm. In the continuous case, on any generic normed spaces $(V, ||\cdot ||_V), (W, ||\cdot ||_W) $, we are using  the operator norm $|| A || = sup_{u \neq 0} \dfrac{ || Au ||_W}{||u ||_V}$. In the following analysis and the convergence theorems below, one can think of discrete norms. For any particular model problem, one would have different (then specified) continuous norms. The choice of $\alpha_\AA$ gives
\[
\left| b(v,u) + c(u,v)\right| 
\leq \left(\alpha_a-\alpha_\AA\right)\|u\|^2 + \left(\alpha_d-\alpha_\AA\right)\|v\|^2\,,
 \]
and, therefore,
\[
 b(v,u) + c(u,v)\geq \left(\alpha_\AA-\alpha_a\right)\|u\|^2 + \left(\alpha_\AA-\alpha_d\right)\|v\|^2\,,
\]
 where this last inequality implies the coercivity of $\AA$.

\section{Block-iterative methods for coupled PDEs}
\label{sec:iterative_methods}

Classic static iterative methods for linear systems can be adapted to block matrices.
In our case, the block pattern of the coupled operator $\mathcal{A}$ in \cref{eq:system} is determined by the number of equations and unknowns, while the dimensions and properties of each block are determined by the mesh size and discretisation scheme adopted for that particular equation.
Convergence analysis for the particular case of positive definite matrices is treated in \citep[p.~47]{Hackbusch2016}. The same author provides rates of convergence in the case of the Poisson equation. In the following paragraphs, we recall the classical (static) iterative methods in their block extensions. The superscript index indicates the iterative step. In general, any iterative splitting/decoupling scheme can be represented by a decomposition of the operator $\mathcal{A}$ into an implicit part $\Ai$ and an explicit part $\Ae$:
\begin{equation}
\label{eq:monolithic_exact}
\Ai w+\Ae w=\mathcal{A}w=\ff.
\end{equation}
This can be rewritten as a stationary iterative scheme as
\begin{equation}
\label{eq:monolithic_iteration}
\Ai w^{k+1}+\Ae w^{k} = \ff.
\end{equation}
While this splitting is often introduced also in monolithic schemes to, for example, derive operator-based preconditioners \cite{axelsson2001survey}, here, for the resulting system to be "decoupled", we require $\Ai$ is a block-triangular operator. For clarity, we now introduce the block version of the most straightforward stationary iterative schemes: Jacobi, Gauss-Seidel, and SOR.
\paragraph{Block-Jacobi}
Block-Jacobi is the simplest iterative method, and its block-version can be formulated as follows:
\begin{align}
	A\uu^{k+1} = \ff_1 - B\vv^{k}, \label{eq:LP.N2.0.Jacobi}\\
	D\vv^{k+1} = \ff_2 - C\uu^{k}, \label{eq:LP.N2.1.Jacobi}
\end{align}
with 
\[
\Ai=
\begin{bmatrix}
	A & 0 \\
	0 & D 
\end{bmatrix}
\,,\qquad
\Ae=
\begin{bmatrix}
	0 & B \\
	C & 0 
\end{bmatrix}.
\]

\paragraph{Block-Gauss-Seidel}
Analogously to the classical Gauss-Seidel method, the block-Gauss-Seidel method (starting from the first equation) reads
\begin{align}
	A\uu^{k+1} = \ff_1 - B\vv^{k},   \label{eq:LP.N2.0.GS}\\
	D\vv^{k+1} = \ff_2 - C\uu^{k+1}, \label{eq:LP.N2.1.GS}
\end{align}
$$
    \Ai=
    \begin{bmatrix}
    	A & 0 \\
    	C & D 
    \end{bmatrix}
    \,,\qquad
    \Ae=
    \begin{bmatrix}
    	0 & B \\
    	0 & 0 
    \end{bmatrix}.
$$

This standard approach is implemented in commercial and applied research codes for solving multi-physics problems.

\paragraph{Block-SOR}
The block version of the SOR method (starting from the first equation) can be written as follows:
\begin{align}
	A\uu^{k+1} = (1-\omega)A\uu^{k} + \omega\left[\ff_1 - B\vv^{k}\right], \label{eq:LP.N2.0.SOR}\\
	D\vv^{k+1} = (1 - \omega)D\vv^{k} + \omega\left[\ff_2 - C\uu^{k+1}\right],\label{eq:LP.N2.1.SOR}
\end{align}
or, equivalently,
\begin{align}
	A\uu^{k+1}
 +  B\vv^{k}
 + \frac{(1-\omega)}{\omega}A(\uu^{k+1}-\uu^{k})
 =
 \ff_1 , \label{eq:LP.N2.0.SOR2}\\
	D\vv^{k+1} 
+
C\uu^{k+1}
+ \frac{(1-\omega)}{\omega}D(\vv^{k+1}-\vv^{k})
 =
 \ff_2 , \label{eq:LP.N2.1.SOR2}
\end{align}
\[
\Ai=
\begin{bmatrix}
	\left(\frac{1}{\omega}\right)A & 0 \\
	C & 	\left(\frac{1}{\omega}\right)D 
\end{bmatrix}
\,,\qquad
\Ae=
\begin{bmatrix}
	\left(\frac{\omega-1}{\omega}\right)A & B \\
	0 & \left(\frac{\omega-1}{\omega}\right)D 
\end{bmatrix}.
\]

\subsection{Convergence properties}
%
\begin{theorem}
The convergence of the iteration \cref{eq:monolithic_iteration} is guaranteed if
\begin{equation}
\label{eq:monolithic_condition_nostab}
\alpha_\AA>2\|\Ae\|\,,
\end{equation}
or
\begin{equation}
    \label{eq:monolithic_condition_nostab2}
    \alpha_{\Ai}>\|\Ae\|\,,
\end{equation}    
where $\alpha_\AA$ is the coercivity constant of the monolithic operator, $\alpha_{\Ai}$ is the coercivity constant of the implicit part, and $\|\Ae\|$ is the continuity constant of the explicit part of the operator.
\end{theorem}
\begin{proof}
From \cref{eq:monolithic_iteration} we can derive the following error equations for $e^k:=w-w^k$:
\begin{equation}
\label{eq:monolithic_erroreqn_nostab}
    \AA e^{k+1} = \Ae (e^{k+1}-e^k)\,;\qquad \Ai e^{k+1} = -\Ae e^k
\end{equation}
If we multiply by $e^{k+1}$ both equations, use the coercivity and the Cauchy-Schwarz inequality we obtain
the following two bounds
\[
\alpha_\AA\|e^{k+1}\| \leq \|\Ae\|\|e^{k+1}-e^k\|\,;
\qquad 
\alpha_{\Ai}\|e^{k+1}\| \leq \|\Ae\|\|e^k\|\,,
\]
i.e.,
\[
\left(\frac{\alpha_\AA}{\|\Ae\|}-1\right)\|e^{k+1}\| \leq \|e^{k}\|\,;
\qquad
\frac{\alpha_{\Ai}}{\|\Ae\|}\|e^{k+1}\| \leq \|e^{k}\|\,;
\]
The fixed-point iterations converge if the LHS coefficient is larger than one, giving the conditions \cref{eq:monolithic_condition_nostab,eq:monolithic_condition_nostab2}.
\end{proof}

\subsection{Relaxation schemes for iterative splitting}
One might find that the block Gauss-Seidel method \cref{eq:LP.N2.0.GS,eq:LP.N2.1.GS} might not converge (for example, it might not satisfy
\cref{eq:monolithic_condition_nostab}). A remedy is to stabilise the equations, similar to what is done in the \textit{fixed-stress} \cite{settari1998, mikelicwheeler, kim2011fixedstress, storvik2019optimization, Both20171June} or \textit{undrained} \cite{kim2011undrained} splitting methods for poro-elasticity equations.
 We introduced here a generalised stabilised/relaxed iteration as follows: Given $\vv^{k}, \uu^{k}$ find $\vv^{k+1}, \uu^{k+1}$ such that
\begin{align}\label{eq:stabilization_0}
    A\uu^{k+1} + B\vv^{k}   + L_u\left(\uu^{k+1}-\uu^{k}\right) &= \ff_1 \\
	C\uu^{k+1} + D\vv^{k+1} + L_v\left(\vv^{k+1}-\vv^{k}\right) &= \ff_2\,,
 \label{eq:stabilization_1}
\end{align}
or, equivalently,
\begin{equation}
    \label{eq:relax}
    \left(\Ai+\LL\right) w^{k+1}=\left(\LL-\Ae\right) w^{k} + \ff\,,
\end{equation}
where
$$
    \LL = 
    \begin{bmatrix}
	L_u & 0 \\
	0 & L_v 
    \end{bmatrix}
$$
is the relaxation operator and, contrarily to what is typically done in literature, is not necessarily a diagonal operator.
It is interesting to notice that there is a strong link between this relaxation and the Block-SOR method where, in the SOR method, the relaxation operators are $\frac{1}{\omega}A$ and/or  $\frac{1}{\omega}D$. Although any relaxation operator can be included in the definition of $\Ai$ and $\Ae$, we treat these terms separately to study their influence better.

One of the difficulties with stabilisation is choosing the appropriate operator $\mathcal{L}$, as its choice will significantly influence the scheme's performance. In \cref{sec:schur}, the relation between the stabilised Gauss-Seidel method and an approximate Schur complement approach will be discussed, providing guidelines to choose the stabilisation operator $\mathcal{L}$. Moreover, in \cref{sec:convergence} convergence proofs for several cases of the stabilised block Gauss-Seidel method \cref{eq:stabilization_0,eq:stabilization_1} are provided, which gives a theoretical justification for the method.

\paragraph{$\ell$-scheme}
A simple but effective choice of relaxation for the operators $L_u$ is to be a multiple of the identity and $L_v=0$:
\[
\LL=
\begin{bmatrix}
	\ell\Id & 0 \\
	0 & 0
\end{bmatrix},
\]
with $\ell>0$. This scheme has been used extensively to solve, for example, poromechanics equations, see, e.g. \cite{mikelicwheeler,Both20171June}. When applied to poromechanics, it is often called \textit{fixed-stress} splitting, and it is the standard over-relaxation approach implemented in many commercial codes to increase the diagonal dominance of the first equation and stabilise the iteration.

\section{Approximate Schur-based methods}
\label{sec:schur}
The Schur complement can be interpreted as the Gaussian elimination formula for block matrices. 
One of the Schur complement's most important applications is for the Saddle Point Problem arising, for example, when solving the incompressible Navier-Stokes equations. The well-known Chorin-Temam algorithm and the analogous Yosida method rely on the block-LU decomposition \citep{viguerie2019effective}, which leads to a Schur complement. These two are special cases of the Uzawa algorithm introduced to solve the SPP \citep{golub2003solving}. A similar method has been applied for the case of coupled energy equations \citep{karki2004application}.

To generalise these approaches for generic block systems, we introduce a general iterative algorithm to stabilise the convergence of iterative coupled schemes with approximate Schur complements by mimicking the incomplete block-LU factorisation. Furthermore, in \cref{sec:approx}, we propose simple diagonal approximations to avoid explicitly computing matrix inverses, although more advanced techniques could also be used \citep{filelis2014generic}. The same strategy has been widely applied to
derive optimal preconditioners \citep{axelsson2009preconditioning}.

It is important to notice that if the exact Schur complement were used, the iterative scheme would converge to the exact solution in one iteration. Therefore, we expect that, by using an approximation, we can stabilise the iterations (ensuring their convergence) and reduce the number of iterations needed for convergence. In this sense, we could also consider these methods as accelerators.

\subsection{Approximate Schur complement factorisation}


Let us consider the decomposition
\begin{equation}
\label{eq:LP.pj.decomposition}
	A\uu = \Aapp\uu + (A-\Aapp)\uu 
\end{equation}
where $\Aapp$ is an approximation of $A$ with the property of being computationally easy to invert, for example, the diagonal of $A$.
We introduce the block-Schur matrix factorisation
$$
    S=
    \begin{bmatrix}
    I & 0\\
    -C\Aapp^{-1} & I
    \end{bmatrix}
$$
which, left-multiplied to the LHS and RHS of the original system \cref{eq:system}, gives

\begin{equation}
\label{eq:LP.pj.exact}
    \begin{bmatrix}
    	A & B\\
    	-C\Aapp^{-1}(A-\Aapp) & D-C\Aapp^{-1}B
    	\end{bmatrix}
    	\begin{bmatrix}
    	\uu\\\vv
    	\end{bmatrix}
    	=
    	\begin{bmatrix}
    	\ff_1\\
    	\ff_2 -C\Aapp^{-1}\ff_1
	\end{bmatrix}
\end{equation}
which is an approximate Schur-based decomposition. If we apply the same technique to the first equation, by using the splitting $D\vv = \mathbb{D}\vv + (D-\mathbb{D})\vv$ and
$$
    \begin{bmatrix}
    I & -B\mathbb{D}^{-1}\\
    0 & I
    \end{bmatrix}
$$
instead of $S$, combining the results we obtain
\begin{equation}
\label{eq:LP.pj.exact2}
    \begin{bmatrix}
    	A-B\mathbb{D}^{-1}C & -B\mathbb{D}^{-1}(D-\mathbb{D})\\
    	-C\Aapp^{-1}(A-\Aapp) & D-C\Aapp^{-1}B
    	\end{bmatrix}
    	\begin{bmatrix}
    	\uu\\\vv
    	\end{bmatrix}
    	=
    	\begin{bmatrix}
    	\ff_1 - B\mathbb{D}^{-1}\ff_2\\
    	\ff_2 -C\Aapp^{-1}\ff_1
	\end{bmatrix}
\end{equation}
which we denote as an alternate approximate Schur-based decomposition.\\

These approximate Schur complements do not provide a decoupled system (i.e., a block-triangular matrix); therefore, the next step is to apply a block-iterative method. If we apply the Gauss-Seidel method,  the identity \cref{eq:LP.pj.exact} results in the following iterative scheme:
\begin{align}
    (D-C\Aapp^{-1}B)\vv^{k+1} &= \ff_2 -C\Aapp^{-1}(\ff_1-(A-\Aapp)\uu^{k}),\label{eq:LP.pj.splitA1}\\
    A\uu^{k+1} &= \ff_1-B\vv^{k+1}. \label{eq:LP.pj.splitA0}
\end{align}
%
Using the fact $A\uu^k+B\vv^k=\ff_1$ from the previous iteration, \cref{eq:LP.pj.splitA0} and \cref{eq:LP.pj.splitA1} can be rewritten in the following form:
\begin{align}
\label{eq:LP.schurL}
    C\uu^{k} + D\vv^{k+1} - C\Aapp^{-1}B\left(\vv^{k+1}-\vv^{k}\right) &= \ff_2,\\
    A\uu^{k+1} + B\vv^{k+1} &= \ff_1,
\end{align}
where we can identify the relaxation operator $\LL$ as
\begin{equation*}
    \LL=
    \begin{bmatrix}
    0 & 0 \\
	0 & -C\Aapp^{-1}B
    \end{bmatrix}.
\end{equation*}

In this formulation, it is interesting to notice the role of the approximate Schur complement $C\Aapp^{-1}B$ as a stabilisation/acceleration term.
%
Similarly, from \cref{eq:LP.pj.exact2}, one can obtain the iterative scheme for the alternate approximate Schur-based version, i.e.
\begin{align}
\label{eq:alternate}
    (A-B\mathbb{D}^{-1}C)\uu^{k+1} &= \ff_1 -B\mathbb{D}^{-1}(\ff_2-(D-\mathbb{D})\vv^{k}),
    \\
    (D-C\Aapp^{-1}B)\vv^{k+1} &= \ff_2 -C\Aapp^{-1}(\ff_1-(A-\Aapp)\uu^{k+1}).
\end{align}
As it will be clarified in the next section, it is desirable to reformulate the problem again as a relaxed iteration, i.e.,
\begin{align}
\label{eq:alternateL}
    B\vv^{k} + A\uu^{k+1} - B\mathbb{D}^{-1}C\left(\uu^{k+1}-\uu^{k}\right) &= \ff_1,\\
    C\uu^{k+1} + D\vv^{k+1} - C\Aapp^{-1}B\left(\vv^{k+1}-\vv^{k}\right) &= \ff_2,
\end{align}
with a relaxation operator
\[
\LL=
\begin{bmatrix}
	-B\mathbb{D}^{-1}C & 0 \\
	0 & 	-C\Aapp^{-1}B
\end{bmatrix}.
\]

This alternate algorithm \cref{eq:alternateL}, however,  written in the relaxation form,  is no longer exactly equivalent to the form directly obtained from the Schur factorisation \cref{eq:alternate}.
Here, we can no longer guarantee that the system is consistent in the coupled sense for each iteration $k$, i.e.,
$A\uu^k+B\vv^k\neq\ff_1$ and $C\uu^k+D\vv^k\neq\ff_2$. In the following, we will focus on the second (relaxation-like) form, although our numerical results showed that the two have negligible differences in the convergence properties.

It is essential to notice that the relaxation operator $\LL$ constructed as above is not always coercive. As discussed in \cref{sec:convergence}, this means that $\LL$ can act as a stabilisation (e.g., in the Quad-Laplacian problem) or as an acceleration (e.g., in the dual-porosity problem).

\subsection{Diagonal approximations of Schur operators}
\label{sec:approx}
In \cref{sec:schur}, we have introduced a general method to solve coupled problems iteratively based on approximated Schur complements. In \cref{sec:convergence}, we have shown that provided some key assumptions for this relaxing operator, the splitting scheme converges. In this section, we introduce practical strategies to approximate the Schur complement, suitable to be easily implemented in PDE toolboxes. These approaches are tested in \cref{sec:results}.
\subsubsection{Schur-based Partial Jacobi}
The first case we consider here is when the approximate Schur-based iteration \cref{eq:LP.schurL} (single) or \cref{eq:alternate} (alternate) is build approximating $A$  by its diagonal (Partial Jacobi). We denote this approach as Schur-based Partial Jacobi (SPJ). In its "alternate" version (applying the relaxation to both equations), the relaxation operator is:
$$
    \LL=
    \begin{bmatrix}
    	-B\mbox{diag}(D)^{-1}C & 0 \\
    	0 & -C\mbox{diag}(A)^{-1}B
    \end{bmatrix}.
$$
The resulting operator is still the product of two potentially full operators and a diagonal. While discretising this operator (or assembling it from existing sparse matrices) in the general case could be not straightforward, the case of the dual-porosity model discussed in \cref{sec:model}, where the operators $B$ and $C$ are diagonal, is particularly simple.

The accuracy of the diagonal approximation is directly linked to the condition number of the matrix and, therefore, to the mesh size. For example, for a Laplacian operator on a uniform mesh, the diagonal approximation has a degenerate spectrum with a single repeated eigenvalue, while the full operator spectrum contains separate eigenvalues (with the ratio between the largest and the smallest eigenvalue being the condition number). This effect will be investigated numerically in \cref{sec:results}.

\subsubsection{Schur-based Double Partial Jacobi}
The SPJ approach could lead to non-sparse matrices, even if $A$, $B$, $C$ and $D$ are sparse. In fact, the products $B\mathbb{D}^{-1}C$ and $C\Aapp^{-1}B$ are not sparse in general. Therefore, we approximate two of the three matrices involved in the two products to keep the sparsity property. By approximating $A$ and $C$, the system reads
$$
    \begin{bmatrix}
        \Aapp & B\\
        \mathbb{C} & D
    \end{bmatrix}
    \begin{bmatrix}
        \uu\\
        \vv
    \end{bmatrix}
    =
    \begin{bmatrix}
        \ff_1-(A-\Aapp)\uu\\
        \ff_2-(C-\mathbb{C})\uu
    \end{bmatrix}
$$
that can then be solved iteratively with block-Gauss-Seidel as
\begin{align}
    (D-\mathbb{CA}^{-1}B)\vv^{k+1} &= \ff_2 -\mathbb{CA}^{-1}(\ff_1+(\Aapp-A)\uu^{k}) + (\mathbb{C}-C)\uu^{k},\label{eq:LP.s2pj.splitA1}\\
    A\uu^{k+1} &= \ff_1-B\vv^{k+1}. \label{eq:LP.s2pj.splitA0}
\end{align}
Analogously to the SPJ case, the system can rewritten in the form
\begin{align}
    C\uu^{k} + D\vv^{k+1} - \mathbb{CA}^{-1}B(\vv^{k+1}-\vv^{k}) &= \ff_2,\\
    A\uu^{k+1} + B\vv^{k+1} &= \ff_1
\end{align}
which is equivalent to \cref{eq:LP.schurL} but with a different stabilisation/acceleration term $\mathbb{CA}^{-1}B$.

If both approximations are diagonal, the resulting matrices are all sparse, and we denote this approach as the double SPJ method (S2PJ).
This can be done only if the coupling operators $B,C$ are square matrices/operators.
In its "alternate" version (applying the relaxation to both equations), the S2PJ relaxation operator is:
$$
    \LL=
    \begin{bmatrix}
    	-\mbox{diag}(B)\mbox{diag}(D)^{-1}C & 0 \\
    	0 & -\mbox{diag}(C)\mbox{diag}(A)^{-1}B
    \end{bmatrix}.
$$
This method is readily applicable to any coupled scalar problem as the resulting system does not involve any product of operators.

When the two operators approximated by a diagonal are of the same type (e.g., the Quad-Laplacian problem, see numerical results in \cref{sec:results}), this method has the advantage of being no longer strongly dependent on the condition number of the matrix $A$, as the spectrum of $BD^{-1}$ is, in fact, clustered and well approximated by the diagonals.

\section{Convergence analysis of relaxed splitting schemes}
\label{sec:convergence}
In this section, we study the convergence properties of the generic splitting scheme with relaxation, \cref{eq:relax}. We will first propose a more general result based on the properties of the monolithic operators $\LL$ and $\AA$, following \cite{atkinson2009theoretical,kirby2010functional}. Then we specialise the analysis for the specific cases of interest, namely the block-symmetric case (such as the dual-porosity problem \cref{eq:darcy1}), and the block-skew-symmetric case (such as the quad-Laplacian \cref{eq:quadLaplacian1} with $m_{uv}=-m_{vu}$). In the following, we will make an extensive use of the following
\begin{lemma}\label{lem:lschemelemma}
Let $S$ a self-adjoint linear operator, then the following identity holds
\begin{equation*}
    \langle x,S(x-y) \rangle = \frac{1}{2}\left(\langle x,Sx\rangle + \langle x-y,S(x-y)\rangle - \langle y,Sy\rangle\right)\qquad\forall x,y.
\end{equation*}
\end{lemma}
\begin{proof}
It is sufficient to expand both sides of the equation and use the fact that $\langle Sy,x\rangle = \langle Sx,y\rangle$ since $S=S^\top$.
\end{proof}

\subsection{Convergence of the monolithic iteration}

The generic relaxed splitting scheme can be written as
\begin{equation}
\label{eq:monolithic_numerical}
    \LL w^{k+1} + \Ai w^{k+1} = \ff - \Ae w^k + \LL w^k.
\end{equation}

\begin{theorem}
Under the assumptions of $\LL$ symmetric and coercive, and $\AA$ coercive, the iteration \cref{eq:monolithic_numerical} converges if
\begin{equation}\label{eq:monolithic_convergence}
\frac{\|\Ae\|^2}{2\alpha_\LL} \leq 2\alpha_\AA - \varepsilon_\LL\alpha_\LL
\end{equation} 
where  $\varepsilon_\LL\geq0$ is the constant such that $\|\LL\|=(1+\varepsilon_\LL)\alpha_\LL$.
Furthermore, the optimal convergence rate, obtained choosing $\alpha_\LL^{opt}=\frac{\|\Ae\|^2}{2\alpha_\AA}$ is
\begin{equation}
    \label{eq:monolithic_rate}
    r(\alpha_\LL^{opt}) = \sqrt{(1+\varepsilon_\LL)\frac{\|\Ae\|^2}{\|\Ae\|^2+2\alpha_\AA^2}}.
\end{equation}
\end{theorem}
\begin{proof}
Subtracting the iteration \cref{eq:monolithic_numerical} from the exact equation \eqref{eq:monolithic_exact}, and substituting $\Ai=\AA-\Ae$, the error equation reads:
\begin{equation}
\label{eq:monolithic_erroreqn}
    \LL(e^{k+1}-e^k) + \AA e^{k+1} = \Ae (e^{k+1}-e^k).
\end{equation}
where $e^k:=w-w^k$. The multiplication of \cref{eq:monolithic_erroreqn} by $e^{k+1}$ and the application of \cref{lem:lschemelemma} leads to
$$\langle e^{k+1}-e^k,\LL(e^{k+1}-e^k)\rangle + 2\langle e^{k+1},\AA e^{k+1}\rangle + \langle e^{k+1},\LL e^{k+1}\rangle = \langle e^k,\LL e^k\rangle + 2\langle e^{k+1},\Ae (e^{k+1}-e^k)\rangle.$$
Thanks to the coercivity of $\AA$ and $\LL$ applied to the LHS of the equation and the Cauchy-Schwarz/$\delta$-Cauchy-Schwarz inequality on the last term of the RHS we obtain
$$\alpha_\LL\|e^{k+1}-e^k\|^2 + 2\alpha_\AA\|e^{k+1}\|^2 + \langle e^{k+1},\LL e^{k+1}\rangle \leq \langle e^k,\LL e^k\rangle + \delta\|e^{k+1}\|^2 + \frac{\|\Ae\|^2}{2\delta}\|e^{k+1}-e^k\|^2,$$
we can collect the terms and obtain
\begin{equation}
\label{eq:monolithic_mainestimate}
    \left(\alpha_\LL-\frac{\|\Ae\|^2}{2\delta}\right)\|e^{k+1}-e^k\|^2 + \left(2\alpha_\AA-\delta\right)\|e^{k+1}\|^2 + \langle e^{k+1},\LL e^{k+1}\rangle \leq \langle e^k,\LL e^k\rangle.
\end{equation}

To make the first term in the LHS of \cref{eq:monolithic_mainestimate} non-negative, we require that
\begin{equation}
\label{eq:monolithic_requirement_1}
   \delta \geq  \frac{\|\Ae\|^2}{2\alpha_\LL} 
\end{equation}
and remove the first term from the LHS, leading to
$$\left(2\alpha_\AA-\delta\right)\|e^{k+1}\|^2 + \langle e^{k+1},\LL e^{k+1}\rangle \leq \langle e^k,\LL e^k\rangle\,.$$
To obtain a bound for the convergence rate, we can use the coercivity and continuity of $\LL$ to obtain
$$\left(2\alpha_\AA + \alpha_\LL - \delta\right)\|e^{k+1}\|^2 \leq \|\LL\|\|e^k\|^2.$$
Since $\alpha_\LL\leq\|\LL\|$, there exists $\varepsilon_\LL\geq0$ such that $\|\LL\|=(1+\varepsilon_\LL)\alpha_\LL$, leading to the following expression
$$\left(2\alpha_\AA + \alpha_\LL - \delta\right)\|e^{k+1}\|^2 \leq (1+\varepsilon_\LL)\alpha_\LL\|e^k\|^2.$$
We can rewrite this expression to emphasise the rate of convergence as follows:
\begin{equation}
    \label{eq:monolithic_rateofconvergence}
    \|e^{k+1}\|\leq\sqrt{\frac{(1+\varepsilon_\LL)\alpha_\LL}{2\alpha_\AA + \alpha_\LL - \delta}}\|e^k\|\,.
\end{equation}
The convergence is guaranteed if
$$0 < \frac{(1+\varepsilon_\LL)\alpha_\LL}{2\alpha_\AA + \alpha_\LL - \delta} < 1$$
hence
$\delta < 2\alpha_\AA - \varepsilon_\LL\alpha_\LL $ and $ \delta < 2\alpha_\AA + \alpha_\LL$. In conjunction with \cref{eq:monolithic_requirement_1} we can write
\begin{equation}
    \label{eq:monolithic_requirement_2}
    0 \leq \frac{\|\Ae\|^2}{2\alpha_\LL} \leq \delta \leq 2\alpha_\AA - \varepsilon_\LL\alpha_\LL\,
\end{equation}
which leads to \cref{eq:monolithic_convergence}.
Since we want to minimise the coefficient in front of $\|e^k\|$ we want to pick the smallest possible value of $\delta$, which is $\|\Ae\|^2/2\alpha_\LL$ according to \cref{eq:monolithic_requirement_2}, hence we chose $\alpha_\LL$ in order to minimise 
$$r(\alpha_\LL) := \sqrt{\frac{2(1+\varepsilon_\LL)\alpha_\LL^2}{2\alpha_\LL^2 + 4\alpha_\AA\alpha_\LL - \|\Ae\|^2}}.$$
The study of this single variate optimal problem leads to the optimal contraction rate 
$$\alpha_\LL^{opt}=\frac{\|\Ae\|^2}{2\alpha_\AA}\,,$$
which gives the rate of convergence \cref{eq:monolithic_rate}.
\end{proof}
\begin{remark}
To recover the classical $\ell$-scheme stabilisation, where $\LL=\ell\Id$, we can choose $\varepsilon_\LL=0$ and $\alpha_\LL=\|\LL\|=\ell>0$ and the converges condition \cref{eq:monolithic_requirement_2} becomes:
\begin{equation} \label{eq:monolithic_requirement_3}
    \alpha_\ell\ge \frac{ \|\Ae\|^2}{2\alpha_\AA}\,,
\end{equation}
with the minimiser of the contraction rate being again the equality.
This means that we can always guarantee the convergence of the scheme by choosing appropriately $\alpha_\LL$ for any given block-iterative approach given by $\Ai$ and $\Ae$.
The presence of a more complex relaxation operator, with $\varepsilon_\LL>0$, seems to penalise the convergence rate. This is because, in this analysis, we have not exploited the operator $\LL$ structure but used only its coercivity and continuity. This result, therefore, cannot fully explain the advantages of using more complex relaxation operators such as the Schur-based factorisation. Nevertheless, this allows us to ensure its convergence.
\end{remark}

The stabilisation $\LL$ could also be included in the definition of $\Ai$ and $\Ae$, and its convergence would be covered by Theorem 1. Nevertheless, this would not allow us to apply the theorem above, and for many practical applications, it is helpful to design a (coercive) extra term to stabilise the splitting. In the case of Schur-based operators $\LL$, the coercivity is not always guaranteed (such as in the dual-porosity problem introduced above and further studied below). In a case of non-coercive $\LL$, Theorem 1 could be applied rewriting \cref{eq:monolithic_condition_nostab} as a condition on the norm of $\LL$:
\begin{equation}
    \label{eq:monolithic_condition_noncoercive}
    \|\LL\|\leq\frac{\alpha_\AA}{2}-\|\Ae\|
\end{equation}

\subsection{Convergence for block-Gauss-Seidel iterations with single relaxation}
The proof above has the advantage of being more general, but it also requires conditions for a coercive and symmetric monolithic relaxation operator $\LL$, chosen based on the continuity constant $\Ae$ and the coercivity constant of the monolithic operator $\AA$. In this section, we consider more in detail the case of the Dual-Porosity model (where $B=C$ or $B=C^\top$, and the Quad-Laplacian (when $B=-C^\top$) and we focus on the block-Gauss-Seidel iterations, with a generic relaxation (e.g., the Schur-based explained above) applied on one equation only, i.e., 
\begin{align}
\label{eq:split2}
Cu^{k} + Dv^{k+1} + {L}(v^{k} - v^{k+1}) &= f_2,\\
\label{eq:split1}
Au^{k+1} + Bv^{k+1} &= f_1.
\end{align}
where ${L}$ is, for example, an approximate Schur complement, such as $ L = C \Aapp^{-1} B$.

\subsubsection{Block-skew symmetric case $C = -B^\top$}
We first consider the block-skew symmetric case, which corresponds to the Quad-Laplacian problem \cref{sec:quadLaplacian}
\begin{theorem}
Let $C=-B^\top$, $L$ symmetric and coercive, such that
\begin{equation}\label{eq:assumption}
{\alpha}_{{L}}\geq \frac{\|B\|^4}{{\alpha}_D{\alpha}_A^2},
\end{equation} 
where  ${\alpha}_{{L}}$, ${\alpha}_{D}$ and ${\alpha}_{A}$ are the coercivity constants of ${L}$, $D$, and $A$ respectively. Then the splitting scheme \cref{eq:split2,eq:split1} converges.
\end{theorem}

\begin{proof}
Define the error functions $e_v^k:=v^k-v$ and $e_u^k = u^k -u$, and subtract the exact equation \cref{eq:system} from \cref{eq:split1} and \cref{eq:split2}.
The scalar product against $e_u^{k+1}$ and  $e_v^{k+1}$) leads to:
\begin{align}
\langle D e_v^{k+1},e_v^{k+1} \rangle - \left\langle L \left(e_v^k-e_v^{k+1}\right),e_v^{k+1}\right\rangle -\left\langle B^\top e_u^{k},e_v^{k+1}\right\rangle &= 0,\\
\left\langle Ae_u^{k+1},e_u^{k+1}\right\rangle + \left\langle Be_v^{k+1},e_u^{k+1}\right\rangle &= 0.
\end{align}
Adding the two equations above gives the equality
\begin{align}
\langle De_v^{k+1},e_v^{k+1} \rangle + \left\langle Ae_u^{k+1},e_u^{k+1}\right\rangle + \left\langle {L}\left(e_v^{k+1}-e_v^{k}\right),e_v^{k+1}\right\rangle + \langle B^\top(e_u^{k+1}-e_u^k),e_v^{k+1}\rangle = 0.
\end{align}
By applying \cref{lem:lschemelemma} we obtain
\begin{align}
&\langle De_v^{k+1},e_v^{k+1} \rangle + \left\langle Ae_u^{k+1},e_u^{k+1}\right\rangle +\frac{1}{2} \left\langle {L}e_v^{k+1},e_v^{k+1}\right\rangle+\\
&+\frac{1}{2}\left\langle{L}\left(e_v^{k+1}-e_v^{k}\right),e_v^{k+1}-e_v^{k}\right\rangle-\frac{1}{2}\left\langle{L}e_v^{k},e_v^{k}\right\rangle+\langle B^\top(e_u^{k+1}-e_u^k),e_v^{k+1}\rangle=0. \nonumber
\end{align}
Now, using Cauchy-Schwarz and Young's inequalities on $\langle B^\top(e_u^{k+1}-e_u^k),e_v^{k+1}\rangle$ we obtain
\begin{align}
    &\langle De_v^{k+1},e_v^{k+1} \rangle + \left\langle A\,e_u^{k+1},e_u^{k+1}\right\rangle +\frac{1}{2} \left\langle {L}e_v^{k+1},e_v^{k+1}\right\rangle + \frac{1}{2}\left\langle{L}\left(e_v^{k+1}-e_v^{k}\right),e_v^{k+1}-e_v^{k}\right\rangle\\ &\leq\frac{1}{2}\left\langle{L}e_v^{k},e_v^{k}\right\rangle + \delta\|e_v^{k+1}\|^2 + \frac{1}{4\delta}\left\|B^\top\left(e_u^{k+1}-e_u^k\right)\right\|^2\nonumber\label{eq:firstestimates}
\end{align}
where $\delta>0$ is free to be chosen.
From \cref{eq:split1}, we have 
\begin{align}
\left\langle A(e_u^{k+1}-e_u^{k}),e_u^{k+1}-e_u^{k}\right\rangle = -\left\langle B(e_v^{k+1}-e_v^{k}),e_u^{k+1}-e_u^{k}\right\rangle,
\end{align}
which by the coercivity of $A$ (with constant ${\alpha}_A$) and Cauchy-Schwarz inequality gives
\begin{align}
{\alpha}_A\|(e_u^{k+1}-e_u^{k})\|^2 \leq \|B\|\|e_v^{k+1}-e_v^{k}\|\|e_u^{k+1}-e_u^{k}\|
\end{align}
and thereby
\begin{align}
\|e_u^{k+1}-e_u^{k}\| \leq \frac{\|B\|}{{\alpha}_A}\|e_v^{k+1}-e_v^{k}\|.
\end{align}
Using the coercivity of $A$, $D$ and ${L}$ (with constants ${\alpha}_A$, ${\alpha}_D$ and ${\alpha}_{{L}}$ respectively) we then get from \cref{eq:firstestimates}
\begin{align}
    \nonumber &{\alpha}_D\|e_v^{k+1}\|^2 + {\alpha}_A\|e_u^{k+1}\|^2 +\frac{1}{2} \left\langle {L}e_v^{k+1},e_v^{k+1}\right\rangle + \frac{{\alpha}_{{L}}}{2}\|e_v^{k+1}-e_v^{k}\|^2\\ 
    &\leq \frac{1}{2}\left\langle{L}e_v^{k},e_v^{k}\right\rangle + \delta\|e_v^{k+1}\|^2 + \frac{1}{4\delta}\|B^\top\|^2\frac{\|B\|^2}{{\alpha}_A^2}\|e_v^{k+1}-e_v^{k}\|^2.
\end{align}
Collecting terms, we then have
\begin{align}\label{eq:preconvergenceargument}
    &({\alpha}_D-\delta)\|e_v^{k+1}\|^2 + {\alpha}_A\|e_u^{k+1}\|^2 +\frac{1}{2} \left\langle {L}e_v^{k+1},e_v^{k+1}\right\rangle
    +
    \left(\frac{{\alpha}_{{L}}}{2}-\frac{1}{4\delta}\frac{\|B\|^4}{{\alpha}_A^2}\right)\|e_v^{k+1}-e_v^{k}\|^2
    \leq\frac{1}{2}\left\langle{L}e_v^{k},e_v^{k}\right\rangle.
\end{align}
Choosing $\delta = \frac{{\alpha}_D}{2}$, and using the assumption \cref{eq:assumption} we sum the equation \cref{eq:preconvergenceargument} from $k=0$ to $k=n$ to obtain 
\begin{align}
    &\frac{{\alpha}_D}{2}\sum_{k=0}^n\|e_v^{k+1}\|^2 + {\alpha}_A\sum_{k=0}^n\|e_u^{k+1}\|^2 +\frac{1}{2} \left\langle {L}e_v^{n+1},e_v^{n+1}\right\rangle \leq\frac{1}{2}\left\langle{L}e_v^{0},e_v^{0}\right\rangle,
\end{align}
and it follows that as $n\to \infty$ the norms of the errors converge to zero; $\|e_v^{n+1}\|^2\to 0$, $\|e_u^{n+1}\|^2\to 0$. 
\end{proof}

\begin{remark}\label{rem:contraction}
  Although the proof above does not provide a convergence rate in the inherent norms, we can obtain one in the weighted ${L}$-norm: $\|x\|_{{L}}^2 = \langle {L}x,x\rangle$. To see this, we first realise that $\|x\|_{{L}}$ is a norm due to the coercivity property of ${L}$. Furthermore, we have the bound
  \begin{equation}
      \|x\|_{{L}}^2\leq \|{L}\|\|x\|^2,
  \end{equation}
  and it follows from equation \cref{eq:preconvergenceargument} that
  \begin{align}
    &\frac{{\alpha}_D}{2\|{L}\|}\|e_v^{k+1}\|_{{L}}^2 + {\alpha}_A\|e_u^{k+1}\|^2 +\frac{1}{2}\| e_v^{k+1}\|^2_{{L}} \leq\frac{1}{2}\|e_v^{k}\|_{{L}}^2.\nonumber
\end{align}
Thereby, we have the contraction
\begin{equation}
    \left(\frac{{\alpha}_D}{\|{L}\|}+1\right)\|e_v^{k+1}\|_{{L}}^2 \leq\|e_v^{k}\|_{{L}}^2.
\end{equation}
\end{remark}

\subsubsection{Block-symmetric case $C=B$}
We study here the case of the dual-porosity and other models with a similar structure, namely with block-symmetric square operators with negative coupling terms.
We first study the case of a generic coercive stabilisation and then consider the case arising from the Schur factorisation, which leads to a non-coercive acceleration term $L$.

\begin{theorem}\label{thm:symmetric_coercive}
Under the assumptions of $C=B$, $A+B$ and $D+B$ coercive, $-B$ coercive and bounded, 
and $L$ symmetric and coercive the solution strategy \cref{eq:split2,eq:split1}, i.e., 
\begin{align}
    Au^{k+1}+Bv^{k+1} &= f_1\label{eq:newsymsplit1}\\
    B u^k + Dv^{k+1} + L\left(v^{k+1}-v^k\right) &= f_2\label{eq:newsymsplit2},
\end{align}
converges.
\end{theorem}
\begin{proof}
Subtracting the exact equation \cref{eq:system} from \cref{eq:newsymsplit2}-\cref{eq:newsymsplit1} we obtain the error equations
\begin{align}
    Ae_u^{k+1}+Be_v^{k+1} &= 0\label{eq:newererrsymsplit1}\\
    B e_u^k + De_v^{k+1} + L\left(e_v^{k+1}-e_v^k\right) &= 0\label{eq:newererrsymsplit2}.
\end{align}
Take the inner product of equation \cref{eq:newererrsymsplit1} with $e_u^{k+1}$ and \cref{eq:newererrsymsplit2} with $e_v^{k+1}$, add the resulting equations and add and subtract $\left\langle Be_v^{k+1},e_v^{k+1}\right\rangle$ + $\left\langle B^\top e_u^{k+1},e_v^{k+1}-e_u^{k+1}\right\rangle$
\begin{align}\label{eq:coercive_noncoercive}
   -\left\langle B\left(e_v^{k+1}-e_u^{k+1}\right),e_v^{k+1}-e_u^{k+1}\right\rangle -
    \left\langle B\left(e_u^{k+1}-e_u^k\right),e_v^{k+1}\right\rangle
    \\+ \left\langle (A+B)e_u^{k+1},e_u^{k+1}\right\rangle + \left\langle(D+B)e_v^{k+1},e_v^{k+1}\right\rangle + \left\langle L\left(e_v^{k+1}-e_v^k\right),e_v^{k+1}\right\rangle &= 0. \nonumber
\end{align}
Using \cref{lem:lschemelemma}, the coercivity of $A+B$ and $B+D$, the coercivity and boundedness of $-B$ and the Cauchy-Schwarz inequality we get
\begin{align}
\label{eq:first2inequality}
   \alpha_B\left\|e_v^{k+1}-e_u^{k+1}\right\|^2+ \alpha_{A+B}\left\|e_u^{k+1}\right\|^2 + \alpha_{D+B}\left\|e_v^{k+1}\right\|^2
   +&\\
   \frac{1}{2}\left\langle L\left(e_v^{k+1}-e_v^k\right),e_v^{k+1}-e_v^k\right\rangle+ \frac{1}{2}\left\langle L e_v^{k+1},e_v^{k+1}\right\rangle
   &\leq
   \frac{1}{2}\left\langle L e_v^k,e_v^k\right\rangle +
    \left\|B\right\|\left\|e_u^{k+1}-e_u^k\right\|\left\|e_v^{k+1}\right\|.\nonumber
\end{align}
Subtracting equation \cref{eq:newererrsymsplit1} at iteration $k$ from the same equation at iteration $k+1$ together with the boundedness of $B$ and coercivity of $A$ gives the inequality
\begin{equation*}
    \left\|e_u^{k+1}-e_u^{k}\right\|\leq \frac{\left\|B\right\|}{\alpha_A}\left\|e_v^{k+1}-e_v^{k}\right\|,
\end{equation*}
which together with Young's inequality with constant $\delta$ and the coercivity of $L$ can be applied to \cref{eq:first2inequality} to obtain
\begin{align*}
   \alpha_B\left\|e_v^{k+1}-e_u^{k+1}\right\|^2 + \left(\alpha_{D+B}-\frac{\left\|B\right\|^4}{2\alpha_A^2\delta}\right)\left\|e_v^{k+1}\right\|^2
   +&\\
    \alpha_{A+B}\left\|e_u^{k+1}\right\|^2+ \frac{\alpha_L}{2}\left\|e_v^{k+1}-e_v^k\right\|^2+ \frac{1}{2}\left\langle L e_v^{k+1},e_v^{k+1}\right\rangle
    &\leq
    \frac{1}{2}\left\langle L e_v^k,e_v^k\right\rangle +
    \frac{\delta}{2}\left\|e_v^{k+1}-e_v^k\right\|^2.
\end{align*}
By choosing $\delta = \frac{\|B\|^4}{\alpha_{D+B}\alpha_A^2}$, and $\alpha_L\geq \delta$ the inequality reduces to
\begin{align*}
   \alpha_B\left\|e_v^{k+1}-e_u^{k+1}\right\|^2+ \alpha_{A+B}\left\|e_u^{k+1}\right\|^2 + \frac{\alpha_{D+B}}{2}\left\|e_v^{k+1}\right\|^2+ \frac{1}{2}\left\langle Le_v^{k+1},e_v^{k+1}\right\rangle &\leq \frac{1}{2}\left\langle L e_v^k,e_v^k\right\rangle,
\end{align*}
and by the arguments from \cref{rem:contraction}, the solution strategy converges.
\end{proof}

\begin{remark}
    It is important to notice that, although the conditions on $A+B$ and $D+B$ seem particularly strong, it must be noticed that, in the dual-porosity model $A[\cdot]=-\divg(m_u \nabla[\cdot]) + \beta[\cdot]$ and $B[\cdot]=-\beta[\cdot]$, therefore $A+B$ is simply the Laplacian operator, and similarly for $D+B$.
    For more general problems, these conditions are necessary for the coercivity of the operator $\AA$.
\end{remark}

\begin{theorem}
Under the assumptions of $C=B$, $A+B$ and $D+B$ coercive, $-B$ coercive and bounded, 
and
\begin{equation*}
    \alpha_{D+B}+\alpha_{-B}\geq 2\|L\| +\frac{\|B\|^2}{4(\alpha_{A+B}+\alpha_{-B})}+\frac{\|B\|^2}{\alpha_A},
\end{equation*}
the solution strategy \cref{eq:split2,eq:split1}, i.e., 
\begin{align*}
    Au^{k+1}+Bv^{k+1} &= f_1,\\
    B u^k + Dv^{k+1} + L\left(v^{k+1}-v^k\right) &= f_2,
\end{align*}
converges. The convergence rate is given by \begin{equation}\dfrac{\|L\| +\frac{\|B\|^2}{4(\alpha_{A+B}+\alpha_{-B})}+\frac{\|B\|^2}{\alpha_A}}{\alpha_{D+B}+\alpha_{-B}- \|L\|}.
\end{equation}
\end{theorem}
\begin{proof}
    The proof follows the same lines as the proof of Theorem~\ref{thm:symmetric_coercive} until equation~\eqref{eq:coercive_noncoercive}. From there, we reorganise and apply coercivity and boundedness properties to obtain
    \begin{align}
   &\alpha_{-B}\left\| e_v^{k+1}\right\|^2 +\alpha_{-B}\left\| e_u^{k+1}\right\|^2 
    + \alpha_{A+B}\left\| e_u^{k+1}\right\|^2 + \alpha_{D+B}\left\|e_v^{k+1}\right\|^2  \\&= -\left\langle L\left(e_v^{k+1}-e_v^k\right),e_v^{k+1}\right\rangle -\left\langle B e_v^{k+1},e_u^{k+1}\right\rangle 
    -\left\langle Be_u^k,e_v^{k+1}\right\rangle. \nonumber
\end{align}
By the Cauchy-Schwarz inequality, we obtain
    \begin{align}\label{eq:crude_bounds}
   &(\alpha_{A+B}+\alpha_{-B})\left\| e_u^{k+1}\right\|^2 + (\alpha_{D+B}+\alpha_{-B})\left\|e_v^{k+1}\right\|^2  \\&\leq \|L\|\|e_v^{k+1}\|^2 +\|L\|\|e_v^k\|\|e_v^{k+1}\|+\|B\| \|e_v^{k+1}\|\|e_u^{k+1}\| 
    +\|B\|\|e_u^k\|\|e_v^{k+1}\|. \nonumber
\end{align}
From equation \eqref{eq:newererrsymsplit1} at iteration $k$ we obtain
\begin{equation}
    \|e_v^k\|\leq \frac{\|B\|}{\alpha_A}\|e_v^k\|,
\end{equation}
which inserted in \eqref{eq:crude_bounds} together with Young's inequality gives
   \begin{align}
   &(\alpha_{A+B}+\alpha_{-B})\left\| e_u^{k+1}\right\|^2 + (\alpha_{D+B}+\alpha_{-B})\left\|e_v^{k+1}\right\|^2  \\&\leq\|L\|\|e_v^{k+1}\|^2 +\left(\|L\|+\frac{\|B\|^2}{\alpha_A}\right)\|e_v^k\|\|e_v^{k+1}\| 
    +(\alpha_{A+B}+\alpha_{-B})\|e_u^k\|^2+\frac{\|B\|^2}{4(\alpha_{A+B}+\alpha_{-B})}\|e_v^{k+1}\|^2. \nonumber
\end{align}
Dividing by $\|e_v^{k+1}\|$ and collecting terms yield
 \begin{equation}
   \left(\alpha_{D+B}+\alpha_{-B}-\frac{\|B\|^2}{4(\alpha_{A+B}+\alpha_{-B})}-\|L\|\right)\left\|e_v^{k+1}\right\| \leq\left(\|L\|+\frac{\|B\|^2}{\alpha_A}\right)\|e_v^k\|.  \nonumber
\end{equation}
Hence, we have a contraction provided
\begin{equation*}
    \alpha_{D+B}+\alpha_{-B}\geq 2\|L\| +\frac{\|B\|^2}{4(\alpha_{A+B}+\alpha_{-B})}+\frac{\|B\|^2}{\alpha_A}.
\end{equation*}
\end{proof}
\begin{remark} We point out that when $L$ is non-coercive, it no longer plays the role of stabilisation, and therefore, it is not strictly necessary, i.e. the convergence is also ensured for $L = 0$.
As we will observe in \cref{sec:results}, a properly constructed non-coercive $L$ in the dual-porosity problem will lead to acceleration.
    In the special case that $L=0$ and $\alpha_{-B} = \|B\|$ we get the condition
    \begin{equation*}
    \alpha_{D+B}+\alpha_{-B}\geq \frac{\alpha_{-B}^2}{4(\alpha_{A+B}+\alpha_{-B})}+\frac{\alpha_{-B}^2}{\alpha_A}.
\end{equation*}

\end{remark}
\section{Numerical Tests}
\label{sec:results}
This section shows illustrative numerical examples of the one- and two-dimensional Dual-Porosity and Quad-Laplacian models. We will test and compare the following schemes, which are labelled as follows:
(unrelaxed) Block-Jacobi iterations ($BJ$), 
(unrelaxed) Block-Gauss-Seidel ($BGS$),
Shur-based Partial-Jacobi on $u$ ($SPJu$),
Shur-based Partial-Jacobi on $v$ ($SPJv$),
Shur-based Partial-Jacobi on all equations (alternate, $SPJa$) and 
Shur-based Double-Partial-Jacobi on $u$, $v$ and both equations (respectively $S2PJ_u$, $S2PJ_v$ and $S2PJ_a$).
%
\paragraph{One-dimensional Finite Volumes solver}
To test the algorithms in a controlled and simple setup, we implemented a one-dimensional finite volume solver in Python for which we can control each step of the discretisation. The finite-volume formulation has a second-order accuracy for flux reconstruction and boundary conditions. Ghost nodes are used to implement boundary conditions. We adopt the method of manufactured solutions for two one-dimensional coupled problems described in the following sections. Both problems are parameterised by the value $\beta$ to change the properties of the resulting monolithic system by tuning the coupling terms. The sparse linear systems derived by internal blocks are solved with the \texttt{sparse.spsolve} method available in the \texttt{scipy} library.
\paragraph{OpenFOAM\textsuperscript{\textregistered} solver}
Together with one-dimensional tests, we provide two-dimensional tests in the CFD-oriented Finite-Volumes platform OpenFOAM\textsuperscript{\textregistered}.
The same models adopted for one-dimensional tests are tested in the unit square domain with Dirichlet and Neumann boundary conditions.\\

Both solvers and the corresponding algorithms are available open-source \citep{roberto_nuca_2022_7457786}.

\subsection{Dual-porosity model}
Here, we present experiments for the Dual porosity model for both one and two dimensions.
\label{subsec:dual_porosity_numerical_tests}

\subsubsection{One-dimensional example}
Equations \cref{eq:darcy1} are tested in one dimension over the interval $[0,\pi]$ with Dirichlet boundary conditions in a grid of 128 cells. The model in one-dimensional form with forcing terms reads
$$ \beta\,(u-v) - \frac{d}{dx}\left(m_u\frac{du}{dx}\right) = f_1,\qquad x\in\,(0,\pi),$$
$$ \beta\,(v-u) - \frac{d}{dx}\left(m_v\frac{dv}{dx}\right) = f_2,\qquad x\in\,(0,\pi).$$
Method of manufactured solutions requires to set the solutions $u$ and $v$. In addition, we also choose non-trivial coefficients $m_u$ and $m_v$. $\beta$ represents the so-called transfer coefficient, but here, it is mainly used to study the efficiency of the proposed algorithms for increasingly coupled and ill-conditioned cases. $\beta$, in fact, controls the coupling between the two equations and the dominance of diagonals in the discretised system. For the mentioned functions, we choose
\begin{align*}
    u(x)   &=\sin(2x),\\
    v(x)   &=e^{-2x},\\
    m_u(x) &=10^4\left[1 + \frac{\sin(2x)}{2}\right],\\
    m_v(x) &=1+\frac{\sin(4x)}{2}.
\end{align*}
The values for Dirichlet boundary conditions are derived by evaluating $u$ and $v$ in $x=0$ and $x=\pi$. The forcing terms $f_1$ and $f_2$ are derived by substituting $u$, $v$, $m_u$ and $m_v$ in the model.
%
%
\begin{figure}[htbp]
	\centering
 	\includegraphics[width=.5\textwidth]{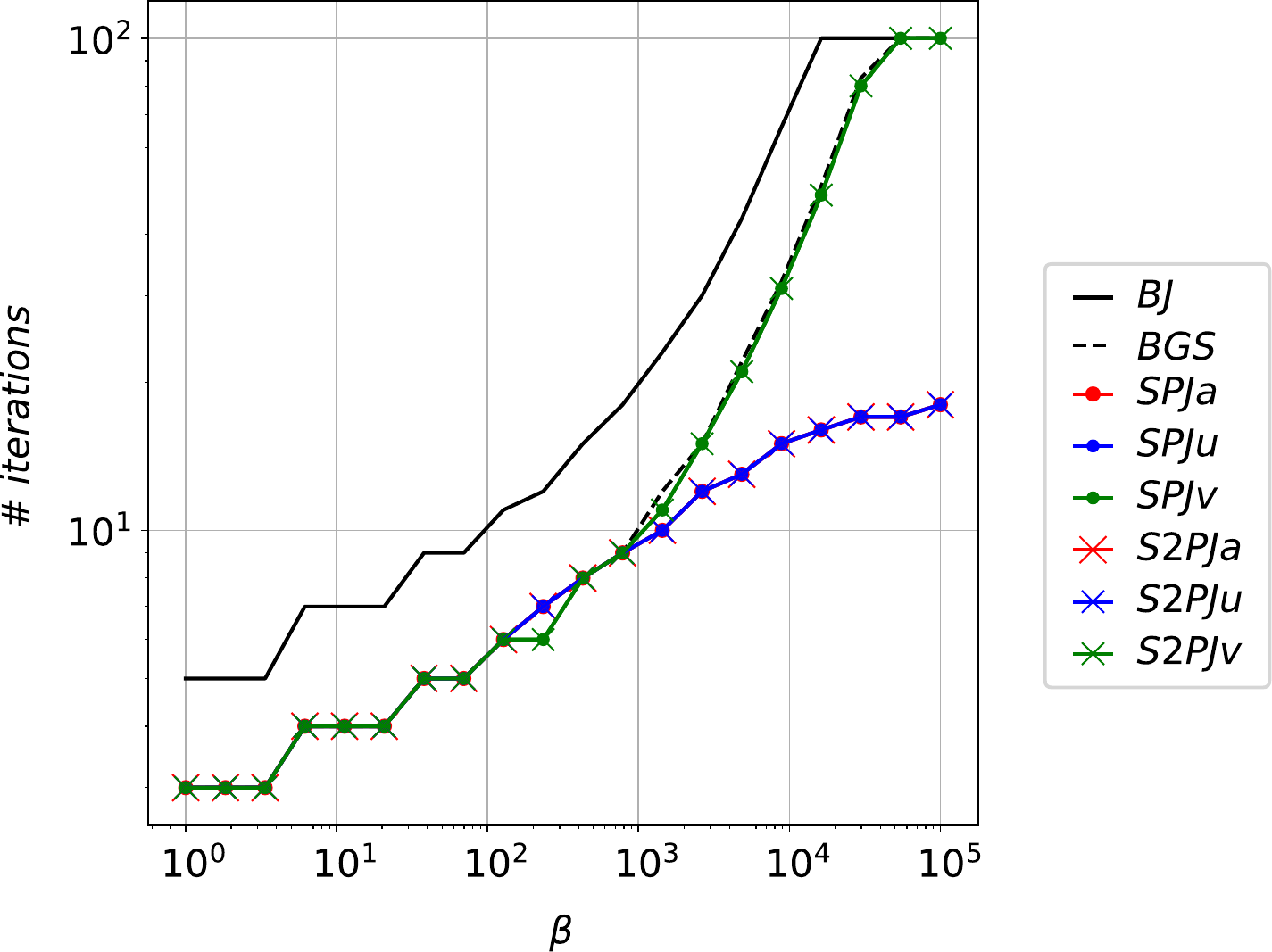}
 	\caption{One-dimensional Dual-porosity problem. Number of iterations required to reach an algebraic residual of $10^{-6}$ for several algorithms as a function of the parameter $\beta$. The plateau at the top of the plot indicates that the maximum number of iterations is reached (100).}
 	\label{fig:dp_1D_iterations}
\end{figure}
In \cref{fig:dp_1D_iterations}, we report the number of iterations needed for the various algorithms to converge up to an algebraic tolerance of $10^{-6}$ or a maximum number of $100$ iterations. As seen, the $SPJ_{(\cdot)}$ class of algorithms behave as well as $BGS$ in the worst-case scenario and much better for $SPJ_a$ and $SPJ_u$. Choosing the equation to apply the relaxation to may be important but not easy. The alternate algorithm has the advantage of being agnostic and applying it to both equations while not compromising its effectiveness. We remark that $SPJ_{(\cdot)}$ and $S2PJ_{(\cdot)}$ coincide because the off diagonal operators are diagonal, hence $\mathbb{B} \equiv B$ or $\mathbb{C} \equiv C$. This case is limited to a very simple mesh to keep it as simple as possible. However, some tests regarding the mesh size are presented for the two-dimensional case.
\FloatBarrier
%
\FloatBarrier
\subsubsection{Two-dimensional example}
\begin{figure}[htbp!]
	\centering
 	\includegraphics[width=.7\textwidth]{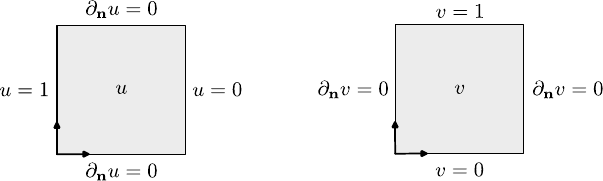}
 	\caption{Two-dimensional Dual-porosity problem. Schematic representations of boundary conditions for $u$ and $v$ in the unit square domain.}
 	\label{fig:dp_2D_boundaryconditions}
\end{figure}

\begin{figure}[htbp!]
	\centering
 	\includegraphics[width=.8\textwidth]{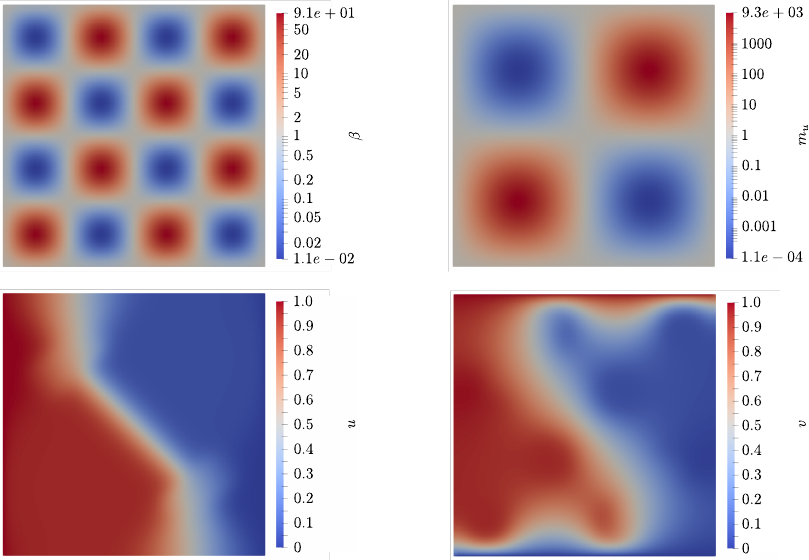}
 	\caption{Two-dimensional Dual-porosity problem. Plots of the function $\beta$, $m_u$ and the solutions $u$ and $v$. We remark that $\beta$ and $m_u$ are plotted in logarithmic scale.}
 	\label{fig:dp_2D_solutions_and_coefficients}
\end{figure}

In this two-dimensional example, we consider two continuous, highly heterogeneous, strongly anti-correlated fields $m_{u}$ and $\beta$. This represents that fracture flow is dominant where the matrix is negligible and vice versa. \Cref{fig:dp_2D_boundaryconditions} shows the boundary conditions for $u$ and $v$ in the square domain. \Cref{fig:dp_2D_solutions_and_coefficients} shows the plots of the function $\beta$, $m_u$ and the solutions $u$ and $v$, the value of $m_v$ is set equal to one and the source terms $f$ are zero.  The mesh is uniform and Cartesian.

\begin{table}[htbp!]
    \centering
    \begin{tabular}{ l | c | c | c}
        method & $50\times50$ & $100\times100$ & $200\times200$ \\
        \hline
        $BJ$     & $90$  & $104$ & $105$\\ 
        $BGS$    & $59$  & $60$  & $57$\\  
        $SPJ_u$ & $43$  & $54$  & $55$\\
        $SPJ_v$ & $31$  & $47$  & $53$\\
        $SPJ_a$ & $28$  & $40$  & $50$
    \end{tabular}
    \caption{Two-dimensional Dual-porosity problem. Number of iterations to reach convergence condition for several algorithms implemented in OpenFOAM\textsuperscript{\textregistered}. Columns indicate the mesh resolutions in $x$ and $y$ directions, respectively, and rows indicate the adopted algorithm.}
    \label{tab:dp_2D_iterations}
\end{table}

\Cref{tab:dp_2D_iterations} shows the number of iterations for each algorithm implemented in OpenFOAM\textsuperscript{\textregistered} for three different meshes. In this simulation, we set $10^{-6}$ as the tolerance for residuals of $u$ and $v$. The results show that $SPJ$ acts as an acceleration term and reduces the number of iterations, with the alternate version showing the best performances. However, this acceleration deteriorates for larger condition numbers of the problem.
%
\FloatBarrier
\subsection{Quad-Laplacian model}
\label{subsec:quad_laplacian_numerical_tests}
This section follows the same testing procedure of the Dual-porosity case for the Quad-Laplacian problem of \cref{eq:quadLaplacian1}. This model has non-trivial off-diagonal blocks, which are discretisations of Laplace operators. This highlights the effect of the approximations of operators in the Schur complement. To our knowledge, there are no direct physical interpretations of this model. However, we can imagine that in a multi-physics context, equations can get easily coupled through diffusive interactions. Therefore, we adopt the Quad-Laplacian as a suitable toy problem for numerical tests. We present one- and two-dimensional tests in the same fashion as the tests presented for the Dual-porosity model. We consider here a Quad-Laplacian problem in one- and two-dimensions with a skew-block-symmetric structure, namely $m_{uv}=-m_{vu}$, in connection with Theorem 3.

\subsubsection{One-dimensional example}
Equations \cref{eq:quadLaplacian1} are tested in one dimension over the interval $[0,2\pi]$ with Dirichlet boundary conditions in a grid of 128 cells.
For the manufactured solutions and model coefficients, we choose
\begin{align*}
    u(x)&=e^{\sin x},\\
    v(x)&=-x^2+x-1,\\
    m_{uu}&=1+\frac{\sin(4x)}{2},\\
    m_{vv}&=\frac{1}{\beta}\left(10^{-2} + \frac{10^{-2}\sin(2x)}{2}\right),\\
    m_{uv}&=-m_{vu}=\beta.
\end{align*}

\begin{figure}[htbp]
	\centering
 	\includegraphics[width=.5\textwidth]{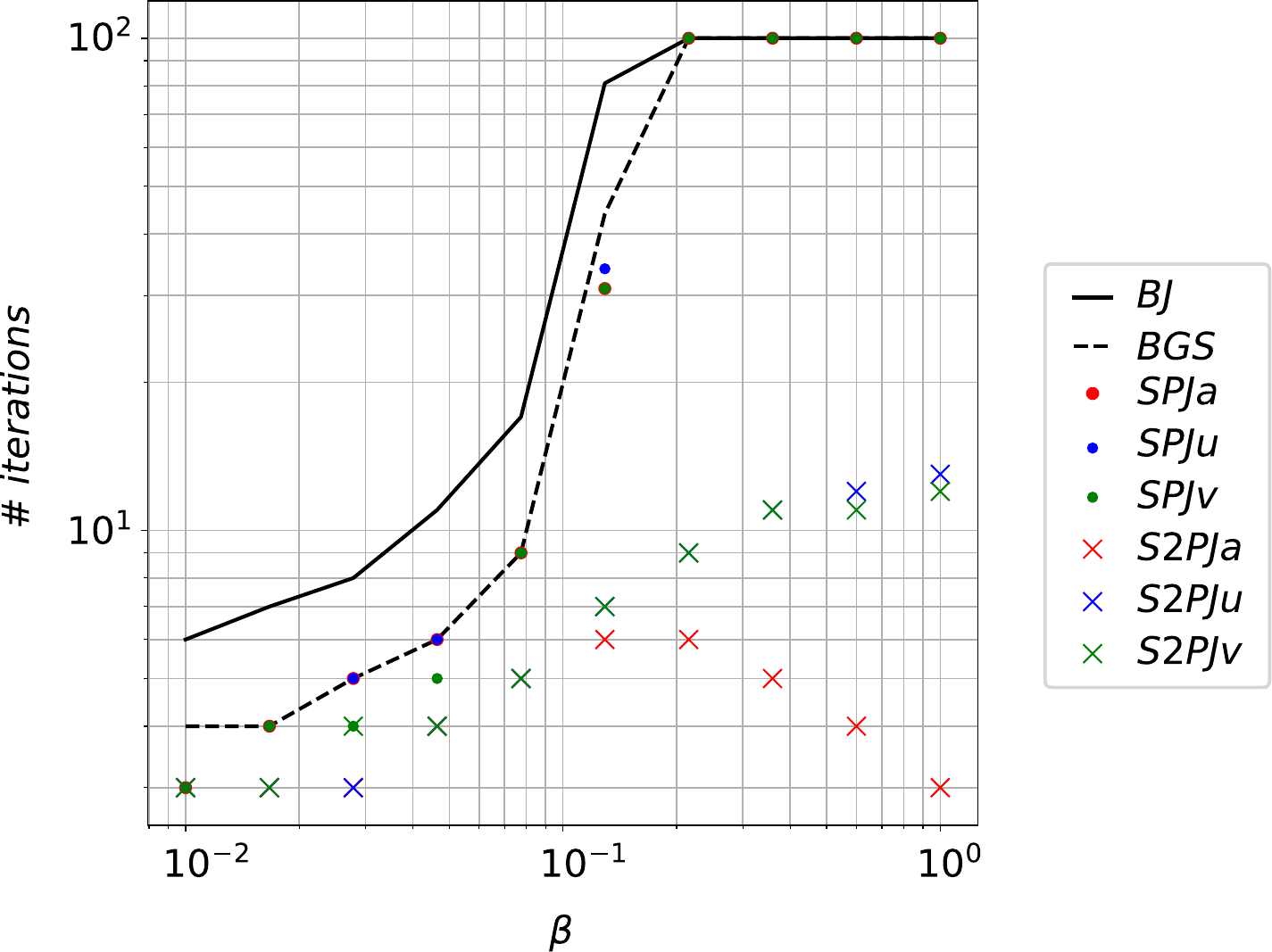}
 	\caption{One-dimensional Quad-Laplacian problem. Number of iterations required to reach an algebraic residual of $10^{-6}$ for several algorithms as a function of the parameter $\beta$. The plateau at the top of the plot indicates that the maximum number of iterations is reached (100).}
 	\label{fig:ql_1D_iterations}
\end{figure}

The values for Dirichlet boundary conditions are derived by evaluating $u$ and $v$ in $x=0$ and $x=\pi$. The forcing terms $f_1$ and $f_2$ are derived by substituting $u$, $v$, $m_{uu}$, $m_{vv}$, $m_{uv}$ and $m_{vu}$ in the model. The parameter $\beta$ appears in the diffusivity coefficients $m_{vv}$, $m_{uv}$ and $m_{vu}$, and it plays the same role as in the Dual-porosity problem. In \cref{fig:ql_1D_iterations}, we report the number of iterations needed for the various algorithms to converge up to an algebraic tolerance of $10^{-6}$ or a maximum number of $100$ iterations. Unlike the Dual-Porosity case, the $S2PJ_{(\cdot)}$ class of schemes behaves much better than $SPJ_{(\cdot)}$.

%
\FloatBarrier
\subsubsection{Two-dimensional example}
\begin{figure}[htbp]
	\centering
 	\includegraphics[width=.8\textwidth]{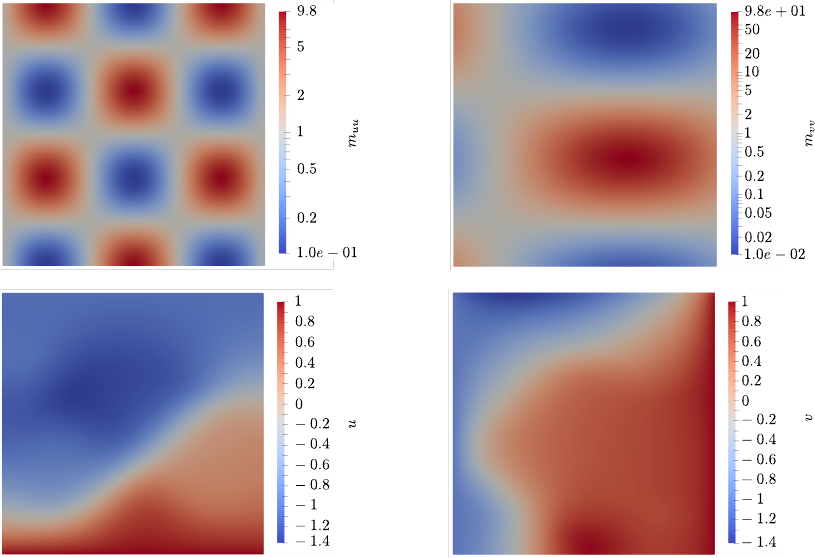}
 	\caption{Two-dimensional Quad-Laplacian problem. Plots of the function $m_{uu}$, $m_{vv}$ and the solutions $u$ and $v$. We remark that $m_{uu}$ and $m_{vv}$ are plotted in logarithmic scale.}
 	\label{fig:ql_2D_solutions_and_coefficients}
\end{figure}

For the two-dimensional case, we have chosen heterogeneous parameters $m_{uu}$ and $m_{vv}$ as shown in \cref{fig:ql_2D_solutions_and_coefficients}, while $m_{uv}=-m_{vu}=1$. The mesh is uniform and Cartesian. OpenFOAM\textsuperscript{\textregistered}

\begin{figure}[htbp]
	\centering
 	\includegraphics[width=.7\textwidth]{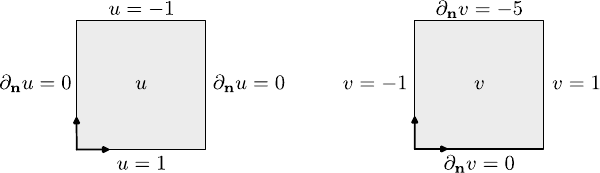}
 	\caption{Two-dimensional Quad-Laplacian problem. Schematic representations of boundary conditions for $u$ and $v$ in the unit square domain.}
 	\label{fig:ql_2D_boundaryconditions}
\end{figure}
\begin{table}[htbp]
    \centering
    \begin{tabular}{ l | c | c | c}
        method & $50\times50$ & $100\times100$ & $200\times200$ \\
        \hline
        $BJ$     & $\times$  & $\times$ & $\times$\\ 
        $BGS$    & $\times$  & $\times$ & $\times$\\  
        $S2PJ_u$ & $37$  & $37$  & $38$\\
        $S2PJ_v$ & $26$  & $26$  & $26$\\
        $S2PJ_a$ & $15$  & $13$  & $13$
    \end{tabular}
    \caption{Two-dimensional Quad-Laplacian problem. Number of iterations to reach the desired residue's tolerance of $10^{-6}$ for several algorithms implemented in OpenFOAM\textsuperscript{\textregistered}.  Columns indicate the mesh resolutions in $x$ and $y$ directions, respectively, and rows indicate the adopted algorithm. The symbol "$\times$" indicates that the algorithm failed to converge.}
    \label{tab:ql_2D_iterations}
\end{table}

Results for different meshes are shown in \cref{fig:ql_2D_solutions_and_coefficients} and in \cref{tab:ql_2D_iterations}.
The method is robust with respect to the mesh size and the conditioning number of the problem. Here, the relaxation operator acts as a stabilisation term, ensuring the convergence of the splitting scheme. Since the high-level interfaces provided by OpenFOAM\textsuperscript{\textregistered} do not allow to compute matrix products of sparse matrices required in the $SPJ$, we can use here the $S2PJ$ variant, thanks to the diagonal approximations of two operators appearing in the double product; only one of the three matrices is sparse, and the other two are simple fields (interpreted as diagonal matrices). However, we also recall that this double approximation makes sense only if the off-diagonal discrete operators are square matrices. 

The S2PJ performs consistently well due to the structure of the problem. Being all operators Laplacian, the two diagonal approximations tend to cancel each other, leading to an accurate approximation of the Schur complement. Therefore, while $BJ$ and $BGS$ diverge, $S2PJ$ converges in a few iterations.
\FloatBarrier
\section{Conclusions}
\label{sec:conclusions}
This work focuses on the development of a unified treatment and theory for iterative splitting schemes for coupled systems of differential equations.
We have shown how stationary iterative methods for linear systems can be applied to systems of Partial Differential Equations (PDEs) and demonstrated the need to introduce a relaxation operator to stabilise the iteration and ensure convergence. We have extended the idea of relaxation to include generic operators and shown how these can
be built based on approximated Schur complements. Convergence theorems are proposed and define sufficient conditions on the relaxation operators.
Numerical tests in one and two dimensions for the Dual-Porosity and Quad-Laplacian problems have been performed using two different open-source codes, one based on the Python library \texttt{scipy} and the other on the OpenFOAM\textsuperscript{\textregistered}) libraries. The numerical results confirm the theory and show the potential of the Schur-based relaxation operator, either as a stabilisation term to ensure convergence or as an acceleration
to reduce the number of iterations, depending on the problem.
The approaches and methods presented here can be extended to include the generalisation to $N\times N$ systems of equations and non-linear problems.
In future works, we aim to include the extension of the convergence theory to better exploit the specific structure of Schur-based relaxation operators.
%
\section{Acknowledgements}
This work has been funded by the following grants:
\begin{itemize}
\item MI has been supported by the European Union’s Horizon 2020 research and innovation programme, grant agreement number 764531, ”SECURe – Subsurface Evaluation of Carbon capture and storage and Unconventional risks”, and by the Nottingham EPSRC Impact Acceleration Award;
\item RN has been supported by the Royal Academy of Engineering, with the grants: "GHG-intensive-industry Medium-scale Capture and Utilization Solutions Assessment In Colombia” and “Asphaltene dynamics at the pore-scale and the impact on oil production at the field-scale”, and by the Nottingham Interdisciplinary Fund.
\end{itemize}

\bibliographystyle{elsarticle-num}
\bibliography{references}

\begin{thebibliography}{10}
\expandafter\ifx\csname url\endcsname\relax
  \def\url#1{\texttt{#1}}\fi
\expandafter\ifx\csname urlprefix\endcsname\relax\def\urlprefix{URL }\fi
\expandafter\ifx\csname href\endcsname\relax
  \def\href#1#2{#2} \def\path#1{#1}\fi

\bibitem{preciceBUNGARTZ2016250}
H.~Bungartz, F.~Lindner, B.~Gatzhammer, M.~Mehl, K.~Scheufele, A.~Shukaev,
  B.~Uekermann, pre{CICE} -- {A} fully parallel library for multi-physics
  surface coupling, Computers \& Fluids 141 (2016).
\newblock \href
  {https://doi.org/https://doi.org/10.1016/j.compfluid.2016.04.003}
  {\path{doi:https://doi.org/10.1016/j.compfluid.2016.04.003}}.

\bibitem{settari1998}
A.~Settari, F.~Mourits, A coupled reservoir and geomechanical simulation
  system, Spe Journal 3~(03) (1998) 219--226.

\bibitem{kim2011fixedstress}
J.~Kim, H.~Tchelepi, R.~Juanes, Stability and convergence of sequential methods
  for coupled flow and geomechanics: {F}ixed-stress and fixed-strain splits,
  Computer Methods in Applied Mechanics and Engineering 200~(13-16) (2011)
  1591--1606.

\bibitem{mikelicwheeler}
A.~Mikeli{\'c}, M.~Wheeler, Convergence of iterative coupling for coupled flow
  and geomechanics, Computational Geosciences 17~(3) (2013) 455--461.

\bibitem{jakubaml}
J.~Both, M.~Borregales, J.~Nordbotten, K.~Kumar, F.~Radu, Robust fixed stress
  splitting for {B}iot's equations in heterogeneous media, Applied Mathematics
  Letters 68 (2017) 101--108.

\bibitem{kim2011undrained}
J.~Kim, H.~Tchelepi, R.~Juanes, Stability and convergence of sequential methods
  for coupled flow and geomechanics: {D}rained and undrained splits, Computer
  Methods in Applied Mechanics and Engineering 200~(23-24) (2011) 2094--2116.

\bibitem{storvik2019optimization}
E.~Storvik, J.~Both, K.~Kumar, J.~Nordbotten, F.~Radu, On the optimization of
  the fixed-stress splitting for {B}iot's equations, International Journal for
  Numerical Methods in Engineering 120~(2) (2019) 179--194.

\bibitem{illiano2021iterative}
D.~Illiano, I.~Pop, F.~Radu, Iterative schemes for surfactant transport in
  porous media, Computational Geosciences 25~(2) (2021) 805--822.

\bibitem{evans2020proof}
C.~Evans, S.~Pollock, L.~Rebholz, M.~Xiao, A proof that {A}nderson acceleration
  improves the convergence rate in linearly converging fixed-point methods (but
  not in those converging quadratically), SIAM Journal on Numerical Analysis
  58~(1) (2020) 788--810.

\bibitem{both2019anderson}
J.~Both, K.~Kumar, J.~Nordbotten, F.~Radu, Anderson accelerated fixed-stress
  splitting schemes for consolidation of unsaturated porous media, Computers \&
  Mathematics with Applications 77~(6) (2019) 1479--1502.

\bibitem{storvik2021accelerated}
E.~Storvik, J.~Both, J.~Sargado, J.~Nordbotten, F.~Radu, An accelerated
  staggered scheme for variational phase-field models of brittle fracture,
  Computer Methods in Applied Mechanics and Engineering 381 (2021) 113822.

\bibitem{douglas1990dual}
J.~Douglas~Jr, T.~Arbogast, Dual porosity models for flow in naturally
  fractured reservoirs, Dynamics of fluids in hierarchical porous media (1990)
  177--221.

\bibitem{fortin1991mixed}
M.~Fortin, F.~Brezzi, Mixed and hybrid finite element methods, Vol.~3, New
  York: Springer-Verlag, 1991.

\bibitem{boffi2013mixed}
D.~Boffi, F.~Brezzi, M.~Fortin, et~al., Mixed finite element methods and
  applications, Vol.~44, Springer, 2013.

\bibitem{hong2021new}
Q.~Hong, J.~Kraus, M.~Lymbery, F.~Philo, A new practical framework for the
  stability analysis of perturbed saddle-point problems and applications,
  Mathematics of Computation (2022).

\bibitem{nicolaides1982existence}
R.~Nicolaides, Existence, uniqueness and approximation for generalized saddle
  point problems, SIAM Journal on Numerical Analysis 19~(2) (1982) 349--357.

\bibitem{atkinson2009theoretical}
W.~Han, K.~Atkinson, Theoretical Numerical Analysis: A Functional Analysis
  Framework, Springer, 2009.

\bibitem{Hackbusch2016}
W.~Hackbusch, {Iterative Solution of Large Sparse Systems of Equations},
  Vol.~95, Springer, 2016.
\newblock \href {https://doi.org/10.1007/978-3-319-28483-5}
  {\path{doi:10.1007/978-3-319-28483-5}}.

\bibitem{axelsson2001survey}
O.~Axelsson, A survey of robust preconditioning methods, in: Topics in
  Numerical Analysis, Springer, 2001, pp. 29--48.

\bibitem{Both20171June}
J.~Both, M.~Borregales, J.~Nordbotten, K.~Kumar, F.~Radu, Robust fixed stress
  splitting for {B}iot's equations in heterogeneous media, Applied Mathematics
  Letters 68 (2017) 101--108.
\newblock \href {https://doi.org/https://doi.org/10.1016/j.aml.2016.12.019}
  {\path{doi:https://doi.org/10.1016/j.aml.2016.12.019}}.

\bibitem{viguerie2019effective}
A.~Viguerie, M.~Xiao, Effective chorin--temam algebraic splitting schemes for
  the steady navier--stokes equations, Numerical Methods for Partial
  Differential Equations 35~(2) (2019) 805--829.

\bibitem{golub2003solving}
G.~Golub, C.~Greif, On solving block-structured indefinite linear systems, SIAM
  Journal on Scientific Computing 24~(6) (2003) 2076--2092.

\bibitem{karki2004application}
K.~Karki, S.~Patankar, Application of the partial elimination algorithm for
  solving the coupled energy equations in porous media, Numerical Heat
  Transfer, Part A: Applications 45~(6) (2004) 539--549.

\bibitem{filelis2014generic}
C.~Filelis-Papadopoulos, G.~Gravvanis, Generic approximate sparse inverse
  matrix techniques, International Journal of Computational Methods 11~(06)
  (2014) 1350084.

\bibitem{axelsson2009preconditioning}
O.~Axelsson, R.~Blaheta, M.~Neytcheva, Preconditioning of boundary value
  problems using elementwise schur complements, SIAM Journal on Matrix Analysis
  and Applications 31~(2) (2009) 767--789.

\bibitem{kirby2010functional}
R.~Kirby, From functional analysis to iterative methods, SIAM review 52~(2)
  (2010) 269--293.

\bibitem{roberto_nuca_2022_7457786}
R.~Nuca, M.~Icardi,
  \href{https://doi.org/10.5281/zenodo.7457786}{splittingschemes v1.0} (2022).
\newblock \href {https://doi.org/10.5281/zenodo.7457786}
  {\path{doi:10.5281/zenodo.7457786}}.
\newline\urlprefix\url{https://doi.org/10.5281/zenodo.7457786}

\end{thebibliography}

\end{document}